\documentclass[11pt,a4paper]{article}
\usepackage{amsmath, amssymb}
\usepackage{latexsym, color}
\usepackage{epsfig}
\usepackage{xcolor}

\numberwithin{equation}{section}

\newtheorem{thm}{Theorem}[section]

\newtheorem{lem}[thm]{Lemma}
\newtheorem{prop}[thm]{Proposition}

\def\nm{\noalign{\medskip}}

\newcommand{\qed}{\hfill \ensuremath{\square}}
\newcommand{\ds}{\displaystyle}
\newcommand{\pf}{\noindent {\sl Proof}. \ }
\newcommand{\p}{\partial}

\newcommand{\norm}[1]{\| #1 \|}

\newcommand{\eqnref}[1]{(\ref {#1})}

\newcommand{\Ibb}{\mathbb{I}}

\newcommand{\Kbb}{\mathbb{K}}
\newcommand{\Rbb}{\mathbb{R}}

\newcommand{\Kcal}{\mathcal{K}}
\newcommand{\Ncal}{\mathcal{N}}

\newcommand{\Scal}{\mathcal{S}}

\def\Ba{{\bf a}}

\def\Bc{{\bf c}}

\def\Bo{{\bf o}}
\def\Bp{{\bf p}}

\def\Bx{{\bf x}}
\def\By{{\bf y}}
\def\Bz{{\bf z}}

\newcommand{\Gd}{\delta}
\newcommand{\Ge}{\epsilon}

\newcommand{\Gvf}{\varphi}

\newcommand{\Gl}{\lambda}

\newcommand{\Gt}{\theta}

\newcommand{\Gs}{\sigma}

\newcommand{\GD}{\Delta}

\newcommand{\GG}{\Gamma}
\newcommand{\GL}{\Lambda}

\newcommand{\beq}{\begin{equation}}
\newcommand{\eeq}{\end{equation}}
\newcommand{\ol}{\overline}

\begin{document}

\title{Quantitative estimates for enhancement of the field excited by an emitter
due to presence of two closely located spherical inclusions\thanks{\footnotesize This work was
supported by NRF grants No. 2015R1D1A1A01059212, 2019R1A2B5B01069967 and 2017R1A4A1014735, by a grant from Central South University and by Hankuk University of Foreign Studies Research Fund of 2019.}}

\author{Hyeonbae Kang\thanks{\footnotesize Department of Mathematics and Institute of Applied Mathematics, Inha University, Incheon 22212, S. Korea, and Department of Mathematics and Statistics, Central South University, Changsha, Hunan, P.R. China (hbkang@inha.ac.kr).} \and KiHyun Yun\thanks{\footnotesize Department of Mathematics, Hankuk University of Foreign Studies, Yongin-si, Gyeonggi-do 17035, S. Korea, and Center for Nonlinear Analysis and Department of Mathematical Sciences, Carnegie Mellon University, Pittsburgh, PA 15213 USA (kihyun.yun@gmail.com).}}

\maketitle

\begin{abstract}
A field in a homogeneous medium can be amplified or enhanced by inserting closely located perfectly conducting inclusions into the medium. In this paper precise quantitative estimates for such enhancement are derived when the given field is the one excited by an emitter of a dipole type and inclusions are spheres of the same radii in three dimensions. Derived estimates reveal the difference, as well as the similarity, between enhancement of the field excited by the emitter and that of the smooth back-ground field. In particular, an estimate shows that when the enhancement occurs, the factor of enhancement is $(\sqrt{\Ge}|\log \Ge|)^{-1}$, which is different from that for the smooth background field, which is known to be $(\Ge|\log \Ge|)^{-1}$ ($\Ge$ is the distance between two inclusions).
\end{abstract}



\section{Introduction and statements of results}

In this paper we investigate enhancement of the field in presence of closely located spherical inclusions in three dimensions, where the field is excited by an emitter located near inclusions. This investigation is motivated by study of the bow-tie shape antenna (see, for example, \cite{PBFLN}) where the electrical field is excited by an emitter, and constitutes a continuation of the earlier work for the cases when inclusions are of bow-tie shape \cite{KY2} and of circular shape in two dimensions \cite{KY3}. Thus we directly go to mathematical description of the problem referring to these papers for motivational remarks and historical accounts.

The problem of this paper can be described in terms of the following mathematical model for $d=2,3$:
\beq\label{main_equation}
\begin{cases}
\GD u = \Ba \cdot \nabla  \Gd_{\Bp}  \quad&\mbox{in } \Rbb^d \setminus \overline {(D_{1} \cup D_{2})}, \\
\ds u =c_j \ (\mbox{constant}) \quad&\mbox{on }\p D_{j},  \\
\ds \int_{\p D_{j} } \p_{\nu} u \, d\Gs  =0, ~&j=1,2,\\
\ds u (\Bx) = O( |\Bx|^{1-d})~&\mbox{as }|\Bx|\rightarrow\infty,
\end{cases}
\eeq
where $D_1$ and $D_2$ are bounded domains representing two inclusions in $\Rbb^d$,  $\nu$ is the unit normal vector on $\p D_1 \cup \p D_2$ pointing {\it inward} to $ D_1 \cup  D_2$, $d\Gs$ denotes the line or surface element on $\p D_j$, and $c_j$ is the constant value attained on $\p D_j$ ($j=1,2$) by the solution $u$. We denote the distance between two inclusions by $\Ge$,
\beq
\Ge = \mbox{dist} (D_1, D_2),
\eeq
and assume that it is small. The term $\Ba \cdot \nabla  \Gd_{\Bp}$ represents the emitter of a dipole type where $\Ba$ and $\Bp$ respectively indicate the direction and the location of the dipole. We assume that $\Bp$ is in the narrow region between two inclusions.
The fact that the solution $u$ attains constant values on $\p D_j$ means that the inclusions are perfect conductors (the conductivity being $\infty$). We emphasize that the constant values $c_j$ are not prescribed but to be determined by the problem. We will show in section \ref{sec:sol} how they are determined.

The solution $u$ to \eqnref{main_equation} exhibits in general singular behavior as $\Ge$ tends to $0$.
To see that, let $\Ncal_\Bp(\Bx):=\GG(\Bx-\Bp)$ where  $\GG$ is the fundamental solution of the Laplacian, {\it i.e.},
\beq\label{gammacond}
\GG (\Bx) =
 \begin{cases}
 \ds \frac{1}{2\pi} \ln |\Bx|\;, \quad & d=2 \;, \\ \nm \ds
 -\frac{1}{4\pi} |\Bx|^{-1}\;, \quad & d = 3 \;.
 \end{cases}
\eeq
Then, $\GD\Ncal_\Bp = \Gd_{\Bp}$, the Dirac delta function, and hence
\beq\label{uzero}
\GD (\Ba \cdot \nabla \Ncal_\Bp) ={\bf a } \cdot \nabla  \Gd_{\Bp}.
\eeq
Thus, $\Ba \cdot \nabla \Ncal_\Bp(\Bx)$ is the solution in absence of inclusions, and its gradient has singularity near $\Bp$ of size $|\Bx-\Bp|^{-3}$ in three dimensions. This singularity in absence of inclusions may be amplified in presence of two inclusions. It is the purpose of this paper to quantify such amplification (or enhancement) of fields when inclusions are spheres of the same radii.

As mentioned earlier, enhancement of fields (excited by an emitter) has been quantitatively studied when the inclusions are of the bow-tie shape with corners separated by the distance $\Ge$ \cite{KY2} and when they are disks \cite{KY3}. It was shown that the singularity of size $|\Bx-\Bp|^{-2}$ in two dimensions is enhanced by the presence of corners (and interaction between them) in the first case, and by interaction of closely located inclusions in the second case. Some of methods to deal with the circular case are adapted in this work to deal with the spherical case. However, there are some significant differences between two and three dimensional cases about which we explain during the course describing results of this paper. Specifically, we remark right after \eqnref{zerozerop} and \eqnref{enhanceorder} on the assumption on the location of the emitter and the order of enhancement, respectively. 

It is insightful to compare the problem of this paper with that of field enhancement when the background field is smooth, particularly, uniform so that it does not have a singularity as the field excited by an emitter does. The problem is described by
\beq\label{main_equation2}
\begin{cases}
\GD u = 0 \quad&\mbox{in } \Rbb^d \setminus \overline {(D_{1} \cup D_{2})}, \\
\ds u =c_j \ \mbox{(constant)} \quad&\mbox{on }\p D_{j},  \\
\ds \int_{\p D_{j} } \p_{\nu} u \, d\Gs  =0, ~&j=1,2,\\
\ds u (\Bx) - h(\Bx) = O( |\Bx|^{1-d})~&\mbox{as }|\Bx|\rightarrow\infty.
\end{cases}
\eeq
Here $h$ is a given harmonic function in $\Rbb^d$ so that $-\nabla h$ represents the background field. The gradient $\nabla u$ of the solution exhibits in general singular behavior in the narrow region between two inclusions as the distance $\Ge$ tends to zero. The problem is to analyse this singular behavior quantitatively. This problem arises in connection with analysis of stress in composites and their effective properties \cite{bab, FK-CPAM-73, Keller-JAP-63}. It has been extensively investigated and several significant results have been obtained in last two decades or so. For examples, it has been shown that the order of blow-up (field enhancement) of $\nabla u$ is $\Ge^{-1/2}$ in two dimensions \cite{AKL-MA-05, Y} and $(\Ge|\log \Ge|)^{-1}$ in three dimensions \cite{BLY-ARMA-09}. We refer to references in a recent paper \cite{KYu} for more comprehensive list of references.

The problems \eqnref{main_equation} and \eqnref{main_equation2}, especially  their unique solvability, can be investigated in a unified way. In fact, we consider the following problem for $f \in H^{1/2}(\p D_1 \cup \p D_2)$ ($H^s$ denotes the Sobolev space of order $s$):
\beq\label{main_equation3}
\begin{cases}
\GD v = 0  \quad&\mbox{in } \Rbb^d \setminus \overline {(D_{1} \cup D_{2})}, \\
\ds v =f- c_j  \quad&\mbox{on }\p D_{j}, \quad j=1,2, \\
\ds \int_{\p D_{j} } \p_{\nu} v \, d\Gs  =0, ~&j=1,2,\\
\ds v (\Bx) = O( |\Bx|^{1-d})~&\mbox{as }|\Bx|\rightarrow\infty
\end{cases}
\eeq
for some constants  $c_j$. We emphasize that this is not an exterior Dirichlet problem with the boundary value $f$. The solution $v$ is in equilibrium in the sense that the total flux across the interfaces is zero as the third condition above states, and constants $c_j$ need to be introduced to achieve such an equilibrium. In fact, the problem is the perfectly conducting problem, namely, the limiting problem (as $k$ tends to $\infty$) of the problem when the conductivity of inclusions, denoted by $k$, is finite. When $k$ is finite, the problem is
\beq\label{finitecond}
\begin{cases}
\GD v = 0  \quad&\mbox{in } (D_{1} \cup D_{2}) \cup (\Rbb^d \setminus \overline {(D_{1} \cup D_{2})}), \\
\ds v|_+ - v|_- = f  \quad&\mbox{on }\p D_{j}, \quad j=1,2, \\
\ds \p_\nu v|_+ - k \p_\nu v|_- = 0  \quad&\mbox{on }\p D_{j}, \quad j=1,2, \\
\ds v (\Bx) = O( |\Bx|^{1-d})~&\mbox{as }|\Bx|\rightarrow\infty.
\end{cases}
\eeq
If we let $v_k$ be the solution to \eqnref{finitecond}, then we see immediately that $\int_{\p D_{j} } \p_{\nu} v_k|_\pm \, d\Gs  =0$ for $j=1,2$, and hence the third condition in \eqnref{main_equation3} follows. If $k \to \infty$, we see formally that $\p_\nu v_k|_- \to 0$, and hence $v_k|_- \to c_j(\mbox{constant})$.

In the next section, we show that the solution to \eqnref{main_equation3} is unique and can be represented in terms of single layer potentials, and as consequences prove unique solvability of problems \eqnref{main_equation} and \eqnref{main_equation2}. Such an approach based on layer potential techniques work quite well for unique solvability. However, in this paper we are interested in the local pointwise behavior of the solution's gradient in the narrow region in between two inclusions. Since the layer potential is not a local operator, it is not clear how the integral equation approach can be utilized for the purpose of this paper. It would be quite interesting to develop effective methods for analysis of field enhancement in the unified context like \eqnref{main_equation3} and the corresponding integral equation. In fact, there has been such an attempt: in \cite{ABTV-ASENS-15} the field enhancement in the context of \eqnref{main_equation2} in two dimensions has been analyzed using an integral equation approach.

In this paper we do not pursue to derive estimates in the general geometry of inclusions. Our purpose is to understand how the enhancement occurs and what the order of the enhancement is when it occurs, in a simple geometry of inclusions. Here we assume that inclusions are balls of the same radii. Furthermore, the fact that inclusions are balls of the same radii is crucially used in the course of deriving estimates. Thanks to the simple geometry of inclusions, we are able to see clearly how the field enhancement differs depending on the direction $\Ba$ of the emitter. In fact, we show that if the direction $\Ba$ is perpendicular to the shortest line between two inclusions (with the location on the perpendicular line), then no enhancement occurs, and if $\Ba$ is parallel to the shortest line, then enhancement actually occurs and the order of enhancement is $(\sqrt{\Ge}|\log \Ge|)^{-1}$. We emphasize that this order is different from that for the problem \eqnref{main_equation2}, which is $(\Ge|\log \Ge|)^{-1}$ as proved in \cite{BLY-ARMA-09}. It is a completely unexpected result to which we address more after stating Theorem \ref{thm100} below.

After translation and rotation if necessary we assume that
\beq\label{config}
D_j = B_{1}((-1)^j (1+\Ge/2),0,0), \quad j=1,2.
\eeq
Here and throughout this paper, $B_r (\Bc)$ denotes the open ball of radius $r$ centered at $\Bc$. If $\Bc$ is the origin $\Bo$, we simply write $B_r$. We assume that the emitter is located on the $x_3$-axis, {\it i.e.},
\beq\label{zerozerop}
\Bp = (0,0, p)
\eeq
for some $p$ with $|p| \leq M$ for some $M$. If $\Bp$ is at some fixed distance from the origin regardless of $\Ge$, then for $\Bx$ in the narrow region between two inclusions, the singularity $|\Bx-\Bp|^{-3}$ which is caused by presence of the emitter is nothing but of order $1$ if $\Ge$ is sufficiently small. Thus it is only meaningful if $M$ becomes smaller as $\Ge$ does. We will assume that $M$ is of size $|\log \Ge|^{-2}$.  We emphasize that this is an improvement over the result in \cite{KY3} for the two-dimensional case where $M$ is assumed to be of size $\sqrt{\Ge}$. We then investigate by quantitative estimates whether the field, the gradient of the solution, is enhanced beyond the singularity $|\Bx-\Bp|^{-3}$.

In fact, the field enhancement takes place because of the potential difference of the solution on the boundaries of the inclusions, namely, $|u|_{\p D_2} - u|_{\p D_1}|$. Moreover, such a potential difference may or may not occur depending on the direction $\Ba$ as mentioned before. Taking advantage of symmetry of the problem, we consider three different cases, namely, those when $\Ba=(1,0,0)$, $\Ba=(0,1,0)$ and $\Ba=(0,0,1)$, and derive different estimates for $|\nabla u|$ for each case. Derived estimates show that when $\Ba$ is either $(0,1,0)$ or $(0,0,1)$, then there is no enhancement of the field, while there is when $\Ba$ is $(1,0,0)$.

If $\Ba=(1,0,0)$, then we see that $u|_{\p D_2} =- u|_{\p D_1} (\neq 0)$, that is, there is actually a potential difference in this case. We obtain different estimates in four regions of interest: when $\Bx$ is close to $\Bp$ ((i) below), $\Bx$ is somewhat away from $\Bp$ ((iii) and (iv) below), and in between ((ii) below). Here and throughout this paper we use the following notation: $A \lesssim B$ means that $A \le C B$ for some constant $C$ independent of $\Ge$ and $\Bp$, and $A \equiv B$ means that both $A \lesssim B$ and $B \lesssim A$ hold.

\begin{thm}\label{thm100}
Suppose that inclusions are of the form \eqnref{config} and let $u$ be the solution to \eqref{main_equation} when $\Ba=(1,0,0)$. There are positive constants $A$, $\Ge_0$, $C_0$ and $C$ such that the following estimates hold for all $\Bp=(0,0,p)$ with $|p| \le C_0|\log \Ge|^{-2}$ and for all $\Ge \le \Ge_0$:
\begin{itemize}
\item[\rm (i)] If $0<|\Bx-\Bp| \leq C (\Ge + |\Bp|^2 )$, then
\beq\label{near}
|\nabla u (\Bx)| \simeq \frac 1 {|\Bx-\Bp|^3}.
\eeq

\item[\rm (ii)] If $C  (\Ge + |\Bp|^2 )   \leq |\Bx-\Bp| \leq C^{-1}  (\Ge + |\Bp|^2 ) |\log \Ge|$, then
\beq\label{between}
|\nabla u (\Bx)| \lesssim \frac 1 {|\Bx-\Bp|^3} \exp \left(-A \frac {|\Bx-\Bp|}{  \Ge  + |\Bp|^2}\right)   + \frac 1 {(\Ge+ |\Bp|^2)^2 |\log \Ge|}.
\eeq

\item[\rm (iii)] If $C^{-1}  (  \Ge + |\Bp|^2 ) |\log \Ge| \leq |\Bx-\Bp| \leq |\log \Ge|^{-2}$, then
\beq\label{far}
|\nabla u (\Bx)| \simeq \frac 1 { (\Ge+|\Bp|^2) |\log \Ge| (\Ge + |\Bx|^2)}.
\eeq

\item[\rm (iv)] If $ |\log \Ge|^{-2} \leq |\Bx-\Bp| $, then
\beq\label{farther}
|\nabla u (\Bx)| \lesssim \frac 1 {(\Ge+|\Bp|^2) |\log \Ge| (|\Bx|^2 + |\Bx|^3)}.
\eeq
\end{itemize}
\end{thm}

Note that \eqnref{near} is an estimate of the singularity produced by the emitter, while \eqnref{between} is a kind of an interpolation of the estimates of the type $|\Bx-\Bp|^{-3}$ and $((\Ge+ |\Bp|^2) |\log \Ge|(\Ge+ |\Bx|^2) )^{-1}$, since $\Ge+ |\Bx|^2 \simeq \Ge+ |\Bp|^2$ in this case. The estimate \eqnref{far} shows that field enhancement actually occurs. For example, if $|p|$ is of order $\sqrt{\Ge}$ and $\Bx$ is on the shortest line segment between two inclusions, namely, $[-\Ge/2, \Ge/2]$ on the $x_1$-axis,  then $|\Bx-\Bp| \simeq \sqrt{\Ge}$, and thus the singularity produced by the emitter is of order $\Ge^{-3/2}$. However, \eqnref{far} shows that
\beq\label{enhanceorder}
|\nabla u (\Bx)| \simeq \frac 1 {\Ge^2 |\log \Ge|}.
\eeq
So, the field is enhanced by the factor of $(\sqrt{\Ge} |\log\Ge|)^{-1}$. It is quite interesting to observe that this factor of enhancement is different from that for the problem \eqnref{main_equation2}, where it is $(\Ge |\log\Ge|)^{-1}$ (see \cite{BLY-ARMA-09}).  In two-dimensional case treated in \cite{KY3}, the factor of enhancement for \eqnref{main_equation} and \eqnref{main_equation2} are the same: it is $\Ge^{-1/2}$. This discrepancy between two and three dimensions may be best described by singular functions, the function $q$ defined in \eqnref{qdef} later: In two dimensions, the singular function (for circular inclusions) is nothing but a sum of two point charges. However, in three dimensions, as we explain in section \ref{sec:sin}, it is infinite sum of point charges. We refer to Lemma \ref{lem:Estimates_for_Q} and its proof for details.

If $\Ba=(0,1,0)$, then one can see from the symmetry of the problem  that $u|_{\p D_2} = u|_{\p D_1}=0$ (see section \ref{sec:decom}). In this case, we obtain the following theorem:

\begin{thm}\label{thm010}
Suppose that inclusions are of the form \eqnref{config} and let $u$ be the solution to \eqref{main_equation} when $\Ba  =  (0,1,0)$. There are positive constants $A$, $\Ge_0$ and $C$  such that the following estimate holds for all $\Bp=(0,0,p)$ with   $|p| \leq  1/2 $ and for all $\Ge \le \Ge_0$:\begin{itemize}
\item[\rm (i)]  If $|\Bx-\Bp|\leq C (\Ge + |\Bp|^2)$, then
\beq\label{est010_2}
|\nabla u (\Bx)| \simeq  \frac 1 {|\Bx - \Bp |^3} .
\eeq
\item[\rm (ii)] If $\Bx \in \Rbb^3 \setminus \overline {D_1 \cup D_2 \cup \{\Bp\}}$, then
\beq\label{est010}
|\nabla u (\Bx)| \lesssim  \frac 1 {|\Bx - \Bp |^3} \exp \left(- \frac {A|\Bx-\Bp|}{ (\sqrt \Ge+|\Bp|)| \Bx- \Bp| + \Ge +|\Bp|^2 }\right).
\eeq
\end{itemize}
\end{thm}

It follows immediately from \eqnref{est010} that
\beq\label{est010-2}
|\nabla u (\Bx)| \lesssim  \frac 1 {|\Bx - \Bp |^3} ,
\eeq
which means that the field is not enhanced beyond $|\Bx - \Bp |^{-3}$. Actually \eqnref{est010} says much more than that. It shows that $\nabla u(\Bx)$ actually decays to $0$ exponentially fast as $\Ge$ tends to $0$ for any $\Bx \neq 0$. In fact, if $\Bx \neq 0$, then $|\Bx-\Bp| \ge \frac{1}{2}|\Bx|$ if $\Ge$ is sufficiently small. Thus $\nabla u(\Bx)$ tends to $0$ exponentially fast as $\Ge$ tends to $0$.

If $\Ba=(0,0,1)$, then $u|_{\p D_2} = u|_{\p D_1}(\neq 0)$, and we obtain the following theorem from which we see as before that the field is not enhanced.

\begin{thm}\label{thm001}
Suppose that inclusions are of the form \eqnref{config} and let $u$ be the solution to \eqref{main_equation} when $\Ba  = (0,0,1)$. There are positive constants $A$, $\Ge_0$ and $C$ such that the following estimate holds for all $\Bp=(0,0,p)$ with   $|p| \leq 1/2$ and for all $\Ge \le \Ge_0$: \begin{itemize}
\item[\rm (i)]    If $|\Bx-\Bp|\leq C (\Ge + |\Bp|^2)$, then \beq\label{est0013} |\nabla u (\Bx)| \simeq  \frac 1 {|\Bx - \Bp |^3} .\eeq
\item[\rm (ii)] If $\Bx \in B_{4} \setminus \overline {D_1 \cup D_2 \cup \{\Bp\}}$, then
\beq\label{est001}
|\nabla u (\Bx)| \lesssim \frac 1 {|\Bx - \Bp|^3} \exp \left(- \frac {A|\Bx- \Bp|}{   (\sqrt \Ge+|\Bp|) |\Bx - \Bp | +   \Ge + |\Bp|^2}\right)+ \exp \left(-\frac{A}{\sqrt \Ge +|\Bx|}\right).
\eeq

\item[\rm (iii)]  If $|\Bx| \ge 4$, then
\beq\label{est0012}
|\nabla u (\Bx)| \lesssim \frac 1 {|\Bx|^3}.
\eeq
\end{itemize}
\end{thm}

\section{An integral equation approach}\label{sec:sol}

In this section, we discuss an integral equation method to address the question of well-posedness of the problem \eqref{main_equation3} (and hence \eqnref{main_equation} and  \eqnref{main_equation2}).

If $D$ is a bounded simply connected domain in $\Rbb^d$ ($d=2,3$) with a Lipschitz boundary, then the single layer potential
$\Scal_{\p D}[\Gvf]$ of a function $\Gvf \in L^2 (\p D)$ is defined by
\beq
\Scal_{\p D} [\Gvf] (\Bx) = \int_{\p D} \GG(\Bx-\By) \Gvf(\By) \, d\Gs (\By), \quad \Bx \in \Rbb^d,
\eeq
where $\GG(\Bx-\By)$ is the fundamental solution to the Laplacian as given in \eqnref{gammacond}.
The following jump formula is well-known:
\beq\label{jump}
-\p_\nu \Scal_{\p D} [\Gvf] \big|_{\pm} (\Bx) =
 \left(\pm {\frac 1 2} I +\Kcal_{\p D}^* \right) [\Gvf]  (\Bx)~~\mbox{a.e. }  \Bx\in \p D,
 \eeq
where
 \beq \label{def_Kcal_D}
 \Kcal_{\p D}^* [\Gvf] (\Bx) = - \int_{\p D} \p_{\nu_\Bx} \GG(\Bx-\By) \Gvf(\By)~d\Gs (\By), \quad \Bx \in \p D.
 \eeq
Here, $\p_\nu$ denotes the normal derivative and the subscripts $+$ and $-$ represent the
limits from outside and inside $D$, respectively. The $-$ signs in the left-hand side  of \eqnref{jump} and  the right-hand one  of \eqnref{def_Kcal_D} are required since the normal vector on $\p D$ is pointing inward in this paper. The operator $\Kcal_{\p D}^*$ is called the Neumann-Poincar\'e (NP) operator.

We seek a solution to \eqnref{main_equation3} in the form of
\beq\label{urepres}
v (\Bx)= \Scal_{\p D_1}[\Gvf_1] (\Bx) +  \Scal_{\p D_2}[\Gvf_2] (\Bx)
\eeq
for a pair of functions $(\Gvf_1, \Gvf_2) \in H^{-1/2}_0(\p
D_1) \times H^{-1/2}_0(\p D_2)$ ($H^{-1/2}_0$ denotes the set of $H^{-1/2}$ functions with the mean zero). If $v$ takes the form \eqnref{urepres}, then all the conditions except the second one in \eqnref{main_equation3} are fulfilled. In fact, the first condition is satisfied since $\Scal_{\p D_j}[\Gvf_j]$ is harmonic in $\Rbb^d \setminus \p D_j$. The last one is fulfilled since $\Gvf_j \in H^{-1/2}_0(\p D_j)$, and hence $\Scal_{\p D_j}[\Gvf_j](\Bx)= O(|\Bx|^{1-d})$. The third condition can be seen to be satisfied using \eqnref{jump} and the fact that $  \Kcal_{\p D_j}[1]=1/2$ on $\p D_j$.

Let us now show that there is a unique pair of potential functions $(\Gvf_1, \Gvf_2)$ such that $u$ given by \eqnref{urepres} is indeed the solution to \eqnref{main_equation3}. For that purpose, let $\GL_j$ be the Dirichlet-to-Neumann map for the domain $D_j$, namely, $\GL_j(f)= \p_\nu u_j|_{\p D_j}$, where $u_j$ is the solution to the Dirichlet problem $\GD u_j=0$ in $D_j$ and $u_j=f|_{\p D_j}$ on $\p D_j$. Then, the second condition in \eqnref{main_equation3} is satisfied if
\beq\label{inteqn1}
\p_\nu (\Scal_{\p D_1}[\Gvf_1] +  \Scal_{\p D_2}[\Gvf_2])|_- = \GL_j(f) \quad\mbox{on } \p D_j, \ j=1,2.
\eeq
In fact, the above relation implies that
$$
\Scal_{\p D_1}[\Gvf_1] +  \Scal_{\p D_2}[\Gvf_2] - f = \mbox{const.} \quad\mbox{on } \p D_j, \ j=1,2.
$$
Thanks to the jump relation \eqnref{jump}, \eqnref{inteqn1} can be rewritten as
\begin{align}
 \left( \frac{1}{2} I - \Kcal_{\p D_1}^* \right) [\Gvf_1] + \p_\nu \Scal_{\p D_2} [\Gvf_2] = \GL_1(f) \quad\mbox{on } \p D_1, \\
\p_\nu \Scal_{\p D_1} [\Gvf_1] + \left( \frac{1}{2} I - \Kcal_{\p D_2}^* \right) [\Gvf_2] = \GL_2(f) \quad\mbox{on } \p D_2,
 \end{align}
or equivalently
\beq\label{Acaleqn}
\left(\frac 12 \Ibb - \Kbb^*\right) \begin{bmatrix} \Gvf_1 \\ \Gvf_2 \end{bmatrix} = \begin{bmatrix} \GL_1(f) \\ \GL_2(f) \end{bmatrix},
\eeq
where
\beq  \mathbb{I}  :=
\begin{bmatrix}  I  & 0 \\
0  & I
\end{bmatrix}
\quad \mbox{and} \quad
\mathbb{K}^*  :=
\begin{bmatrix}  \Kcal_{\p D_1}^* & -\p_\nu \Scal_{\p D_2} \\
-\p_\nu \Scal_{\p D_1}  &  \Kcal_{\p D_2}^*
\end{bmatrix} .
\eeq
It is known that the operator $\frac 1 2\Ibb - \Kbb^*$ is invertible on $H^{-1/2}_0(\p D_1) \times H^{-1/2}_0(\p D_2)$ (see, for example, \cite{ACKLY}). In fact, $\frac{1}{2} I - \Kcal_{\p D_j}^*$ is invertible on $H^{-1/2}_0(\p D_j)$ (see \cite{Verch-JFA-84}) and $\frac 1 2\Ibb - \Kbb^*$ is a compact perturbation of the invertible operator
$$
\begin{bmatrix}  \frac{1}{2} I - \Kcal_{\p D_1}^* & 0 \\
0  & \frac{1}{2} I - \Kcal_{\p D_2}^*
\end{bmatrix}.
$$
Thus invertibility of $\frac 1 2\mathbb{I} - \mathbb{K}^*$ is equivalent to injectivity thanks to the Fredholm alternative. Once \eqnref{Acaleqn} is solved for $(\Gvf_1, \Gvf_2)$, then the function $v$ defined by \eqnref{urepres} yields the unique solution to \eqnref{main_equation3}.

We now show how the constants $u|_{\p D_j}$ can be determined, when $u$ is the solution to \eqref{main_equation2}.  For that, let us recall the following result obtained in \cite{ACKLY}.
\begin{lem}\label{eigenvector}
For $i=1,2$, there is a unique solution $w_i$ to
\beq\label{def_of_vi}
\begin{cases}
\ds \GD w_i  = 0 \quad & \mbox{in } \Rbb^d \setminus \overline{D_1 \cup D_2},\\
\ds w_i = \mbox{constant} \quad &\mbox{on } \p D_j, \ j=1,2,  \\
\nm
\ds \int_{\p D_j} \p_\nu w_i \, d\Gs = \Gd_{ij}, \quad & i, j = 1,2, \\
\nm
\ds w_i(\Bx) - \Gl_i \Scal_{\p D_i} [1](\Bx) = & \hskip-0.3cm O(|\Bx|^{1-d}) \quad \mbox{as } |\Bx|\rightarrow \infty,
\end{cases}
\eeq
for some nonzero constant $\Gl_i$. Here $\Gd_{ij}$ is the Kronecker delta.
\end{lem}

The solution $w_i$ is constructed as $w_i = \lambda_i \Scal_{\p D_i} [1] - v_i$, where $\lambda_i = (\int_{\p D_i} \p_{\nu} \Scal_{\p D_i} [1] |_+ d \sigma)^{-1} $ and $v_i $ is the solution to \eqref{main_equation3} with $f=\lambda_i \Scal_{\p D_i} [1]$.  It is worth mentioning that the solution $w_i$ is closely related to the eigenfunctions of the operator  $\mathbb{K}^*$ corresponding to the eigenvalue $1/2$ ($1/2$ has multiplicity $2$) (see above-mentioned paper).

Having functions $w_1$ and $w_2$ in hand, we define $q$ and $q^\perp$ by
\beq\label{qdef}
q :=-w_1 + w_2 \quad\mbox{and}\quad q^\perp :=w_1 + w_2.
\eeq
Then, we have
\beq\label{int12}
\int_{\p D_j} \p_\nu q \, d\Gs = (-1)^j, \quad\mbox{and}\quad \int_{\p D_j} \p_\nu q^\perp \, d\Gs = 1, \quad j=1,2.
\eeq
Since $w_i$ attains constant values on $\p D_j$, so do $q$ and $q^\perp$. Since $\int_{\p D_1} \p_\nu q \, d\Gs + \int_{\p D_2} \p_\nu q \, d\Gs =0$,
\beq
q(\Bx) =O(|\Bx|^{1-d}) \quad\mbox{as } |\Bx|\rightarrow \infty,
\eeq
while
\beq
q^\perp(\Bx) = \Gl_1 \Scal_{\p D_1} [1](\Bx) + \Gl_2 \Scal_{\p D_2} [1](\Bx) + O(|\Bx|^{1-d}).
\eeq
In particular, $q^\perp(\Bx) = O(|\Bx|^{-1})$ in three dimensions. The function $q$ plays a crucial role for analysis of this paper, and so we review some of its important properties in section \ref{sec:sin}.

If $v$ be the solution to \eqnref{main_equation3} when $f= \Ba \cdot \nabla \Ncal_\Bp$ on $\p D_1 \cup \p D_2$, then, $u=\Ba \cdot \nabla \Ncal_\Bp-v$ is the solution to \eqnref{main_equation}. If $v$ is the solution to \eqnref{main_equation3} when $f= h$ on $\p D_1 \cup \p D_2$ where $h$ is the harmonic function in \eqnref{main_equation2}, then $u=h-v$ is the solution to \eqnref{main_equation2}.
In either case, thanks to \eqnref{int12}, we have
$$
u|_{\p D_1} - u|_{\p D_2} =  \int_{\p(D_1 \cup D_2)} (v-f) \p_\nu q  \, d\Gs.
$$
Since $\int_{\p D_j} \p_{\nu} v d\Gs=0$ and $q$ is constant on $\p D_j$, we further have
$$
u|_{\p D_1} - u|_{\p D_2} = \int_{\p(D_1 \cup D_2)} \left[ (v-f) \p_\nu q - q \p_\nu v \right] \, d\Gs .
$$
We then infer using Green's formula that
\beq\label{difference}
u|_{\p D_1} - u|_{\p D_2}  = - \int_{\p(D_1 \cup D_2)} f \p_\nu q  \, d\Gs .
\eeq
Similarly, we have
\beq\label{sum}
u|_{\p D_1} + u|_{\p D_2}  =  \int_{\p(D_1 \cup D_2)} f \p_\nu q^\perp  \, d\Gs .
\eeq
In this way we determine the constants $u|_{\p D_j}$. In fact, $u|_{\p D_j}$ can be expressed in terms of the function $w_j$, but we write in the forms of \eqnref{difference} and \eqnref{sum} for later use.

\section{Construction of $q$ and $q^\perp$ for spherical inclusions}\label{sec:sin}

The function $q$ defined by \eqnref{qdef} is a building bock in describing the behavior of the solution to problem \eqnref{main_equation2}, It was first used in \cite{Y} for investigation of field enhancement (see also \cite{KLY-SIAP-14, LY}). We review important properties of this function in this section since it is also used in an essential way in this paper. We also show how the function $q^\perp$ is constructed.

When $D_1$ and $D_2$ are unit balls, the function $q$ is constructed explicitly as a weighted sum of the differences of the point charges. To depict it, let $R_j$ be the inversion with respect to $\p D_j$ for $j=1,2$, \textit{i.e.},
$$
R_j(\Bx)=\frac{\Bx-\Bc_j}{|\Bx-\Bc_j|^2}+\Bc_j
$$
($\Bc_j$ is the center of $D_j$). Define points $\Bp_k$ for $k=0,1, \ldots$ by
\beq\label{repeat}
\begin{cases}
\ds\Bp_{2k}=(R_2R_1)^{k}\Bc_2, \\
\ds \Bp_{2k+1}= R_2(R_1R_2)^{k}\Bc_1.
\end{cases}
\eeq
Then $\Bp_k$ are contained in $D_2$ and converge monotonically to $\Bp_{\infty}$ as $k \to \infty$, where $\Bp_{\infty}$ is the fixed point of the combined inversion $R_2R_1$. The points $\Bp_k$ are locations of charges.
To define the weight at each $\Bp_k$, let
\beq
\mu_{k}=\frac{1}{|\Bc_1-\Bp_k|}, \quad k=1,2,\ldots,
\eeq
and
\beq
q_{0}=1 \quad \mbox{and} \quad q_{n}= \prod_{k=1}^n\mu_{k},\ n\geq 1 .
\eeq
It is proved in \cite{LY} that the function $q$ when $D_1$ and $D_2$ are balls is given by
\beq\label{h:sameradii}
q(\Bx)=\frac{1}{\sum_{n=0}^{\infty} q_n} \sum_{n=0}^{\infty} q_n \left(- \GG(\Bx-\Bp_n) + \GG(\Bx+\Bp_n) \right).
\eeq

Let $\Bp_n= (p_n,0,0)$. The following properties are derived from the results in \cite{KLY-SIAP-14}. There exists a constant $C>1$ such that
\beq\label{3700}
\left|p_n - \frac 1 {n+1}\right| + \left|q_n - \frac 1 {n+1}\right| \leq C \sqrt \Ge \leq \frac 1 {2 (n+1)}
\eeq
for all $n \leq 1/(C \sqrt{\Ge})$ ((3.9) in that paper),
\beq
p_n \simeq \sqrt \Ge
\eeq
for all $n \geq 1/(C\sqrt{\Ge})$ (Lemma 3.2 (ii), (iii) and (3.9) in that paper), and
\beq\label{3900}
\sum_{n \geq 1/(C\sqrt{\Ge} )} q_n  \lesssim 1
\eeq
((3.11) there). As a consequence, it was proved that
\beq
\sum_{n= 0 }^{\infty}  q_n  \simeq |\log \Ge|. \label{eq:sum_q_n}
\eeq

Now we show how to construct $q^\perp$. Let $v_f$ be the solution to \eqnref{main_equation3} with $f=2\Ncal_{\Bc_1}$ ($\Bc_1$ is the center of $D_1$) on $\p D_1$ and $\p D_2$. Then $q^\perp$ is given by
\beq\label{qperp}
q^\perp = v_f - 2 \Ncal_{\Bc_1} + q.
\eeq
In fact, $q^\perp$ attains constant values on $\p D_j$ since both $v_f - 2 \Ncal_{\Bc_1}$ and $q$ do. And
$$
\int_{\p D_1} \p_\nu q^\perp \, d\Gs = \int_{\p D_1} \left[ \p_\nu v_f - 2 \p_\nu \Ncal_{\Bc_1} + \p_\nu q \right] \, d\Gs = 0+2-1=1.
$$
Similarly we see that $\int_{\p D_2} \p_\nu q^\perp \, d\Gs=1$ since $\Bc_1$ is located outside $D_2$.

\section{Decomposition of the solution}\label{sec:decom}

Let $q$ be the function defined by \eqnref{qdef}, and let $u$ be the solution to the main problem \eqnref{main_equation}. Then both $q$ and $u$ attain constant values on $\p D_j$ ($j=1,2$), which we denote by $q|_{\p D_j}$, etc. Because of the symmetry of $D_1 \cup D_2$ with respect to the $x_2 x_3$-plane, $q|_{\p D_1}= - q|_{\p D_2}$. Define the function $Q$ by
\beq\label{Qdef}
Q(\Bx) := \frac{u|_{\p D_2} - u|_{\p D_1}}{2q|_{\p D_2}} q(\Bx).
\eeq
Then, $Q$ satisfies
\beq
Q|_{\p D_1} = - \frac{u|_{\p D_2} - u|_{\p D_1}}{2} \quad\mbox{and}\quad Q|_{\p D_2} = \frac {u|_{\p D_2} - u|_{\p D_1}}{2}.
\eeq
If we define the function $r$ by
\beq\label{r_def}
u = Q+r ,
\eeq
then it is the solution to
\beq \label{eqn_r}
\begin{cases}
\GD r = \Ba \cdot \nabla \Gd_\Bp  \quad&\mbox{in } \Rbb^3 \setminus \overline {(D_{1} \cup D_{2})}, \\
\ds r = c \quad&\mbox{on }\p D_{j} , \  j=1,2,\\
\ds  r(\Bx)= O\left(|\Bx|^{-2}\right)  &\mbox{as } |\Bx| \rightarrow \infty,
\end{cases}
\eeq
where $c$ is the common constant given by
\beq\label{c}
c = \frac{u|_{\p D_1} + u|_{\p D_2}}{2}.
\eeq

The constants $u|_{\p D_j}$ ($j=1,2$) vary depending on $\Ba$ as we see below:
\begin{itemize}
\item[\rm{(i)}] If $\Ba=(1,0,0)$, then
\beq\label{symm100U}
\Ba \cdot \nabla \Gd_\Bp (-x_1,x_2, x_3)= - \Ba \cdot \nabla \Gd_\Bp (x_1,x_2, x_3).
\eeq
Thus, the solution $u$ to \eqnref{main_equation} enjoys the same symmetry, namely,
\beq\label{symm100}
u (-x_1,x_2,x_3) = -u (x_1,x_2,x_3).
\eeq
In particular, we have $u|_{\p D_1}= -u|_{\p D_2}$, and the constant $c$ given in \eqnref{c} is $0$.

\item[\rm{(ii)}] If $\Ba=(0,1,0)$, then
\beq\label{symm010U}
\Ba \cdot \nabla \Gd_\Bp (x_1,- x_2,x_3)  = - \Ba \cdot \nabla \Gd_\Bp (x_1,x_2,x_3),
\eeq
and thus $u$ satisfies
\beq\label{symm010}
u(x_1,- x_2,x_3) = - u (x_1,x_2,x_3).
\eeq
In particular, we have $u|_{\p D_1}= u|_{\p D_2}=0$. Thus, $Q \equiv 0$ and $c=0$.

\item[\rm{(iii)}] If $\Ba=(0,0,1)$, then
\beq\label{symm001U}
\Ba \cdot \nabla \Gd_\Bp (x_1, x_2,x_3)  =  \Ba \cdot \nabla \Gd_\Bp (-x_1,x_2,x_3),
\eeq
and hence $u$ satisfies
\beq\label{symm001}
u(x_1,x_2,x_3) =  u (-x_1,x_2,x_3).
\eeq
In particular, we have $u|_{\p D_1}= u|_{\p D_2}$. Thus, $Q \equiv 0$.
\end{itemize}

In the next section we show that enhancement of the field occurs in the first case where $Q$ is not zero, and does not occur in the other two cases. It means in particular that field enhancement is due to the potential gap $u|_{\p D_2} - u|_{\p D_1}$.

 If there is no potential gap, then  $u \equiv r$. Even if $Q$ is not zero, there is a region where $\nabla r$ dominates $\nabla u$.  To estimate $\nabla r$, we derive a gradient estimate for $r_0$ in the following lemma, which is an analogue of $r$. A similar estimate was obtained for circular inclusions in two dimensions \cite[Lemma 2.1]{KY3}. We  follow closely the argument for the two-dimensional case.

\begin{lem}\label{lem:case123}
For a unit vector $\Ba$, let  $r_0$ be the solution to 
\beq \label{eqn_r_0_def}
\begin{cases}
\GD r_0 = \Ba \cdot \nabla \Gd_\Bp  \quad&\mbox{in } \Rbb^3 \setminus \overline {(D_{1} \cup D_{2})}, \\
\ds r_0 = 0 \quad&\mbox{on }\p D_{j} , \  j=1,2,\\
\ds r_0(\Bx)= O\left(|\Bx|^{-1}\right)  &\mbox{as } |\Bx| \rightarrow \infty.
\end{cases}
\eeq 
There exist positive constants $A$ and $\Ge_0$ such that
\beq\label{case12:gamma}
|\nabla r_0 (\Bx)|  \lesssim \left(\frac 1 {|\Bx - \Bp|^3} + \frac 1 {|\Bx - \Bp|^2}  \right) \exp \left(- \frac {A |\Bx- \Bp|}{( \sqrt \Ge + |\Bp|)  |\Bx - \Bp | + \Ge + |\Bp|^2 }\right)
\eeq
for all $\Bx \in \Rbb^3 \setminus \overline{D_1 \cup D_2}$, for all $\Bp$ with $|\Bp| \le 1/2$, and for all $\Ge \le \Ge_0$.
\end{lem}

\pf
One can see easily that
$$
r_0 := r- \frac{r|_{\p D_1}}{q^{\perp}|_{\p D_1}} q^{\perp}
$$
is the solution to \eqref{eqn_r_0_def}. The uniqueness follows from the maximum principle.

With help of the transformation $\Phi$ defined by
\beq\label{Phi}
\Phi (\By) =  \frac {\By-\Bp} {| {\By-\Bp}|^2 } +  \Bp,
\eeq
we let $r_0^\dag$ be the Kelvin transform of $ r_0$, namely,
$$
 r_0^\dag (\By) = \frac 1 {|\By - \Bp|}  r_0 \left(\Phi (\By)\right).
$$
We then define $ r_0^*$ by
\beq\label{dagsharp}
 r_0^* (\Bz) =  r_0^\dag (\By),
\eeq
where
\beq\label{xyz}
\By= \Ge_* \Bz + \Bp_*, \quad \Ge_* = (\Ge + p^2 + (\Ge^2 /4))^{-1}, \quad \Bp_* = (1-\Ge_*)\Bp.
\eeq

Note that the transformation $\Bz \mapsto  \Phi(\Ge_* \Bz + \Bp_*)$ maps the region $\Rbb^3 \setminus \overline{(D_{1} \cup  D_{2} \cup \{\Bp\})}$ onto itself. Thus, $r_0^*$ is harmonic in $\Rbb^3 \setminus \overline{(D_{1} \cup  D_{2})}$ except at $\Bp$. Since $r_0(\Bx)=O(|\Bx|^{-1})$ as $|\Bx| \to \infty$, we see that the following limits exist and are the same:
\beq\label{limitp}
\lim_{\Bz \to \Bp} r_0^*(\Bz) = \lim_{\By \to \Bp} r_0^\dag(\By) = \lim_{\Bx \to \infty} |\Bx| r_0 (\Bx)  .
\eeq
Therefore, $\Bp$ is a removable singularity for $r_0^*$. Thus $r_0^*$ satisfies
\beq \label{equation_r_2}
\begin{cases}
\GD  r_0^* = 0  \quad&\mbox{in } \Rbb^3 \setminus \overline{(D_{1} \cup  D_{2})}, \\
\ds r_0^* = 0 \quad&\mbox{on }\p (D_{1} \cup  D_{2}).
\end{cases}
\eeq

We now look into the boundary condition of $r_0^*$ at $\infty$. The function $r_0 - \Ba \cdot \nabla \Ncal_\Bp$ is well-defined in a neighborhood of $\Bp$ as a harmonic function so that
\beq\label{518}
(r_0- \Ba \cdot \nabla \Ncal_\Bp)(\Bx)= O(1) \quad\mbox{as } \Bx \to \Bp,
\eeq
or equivalently $(r_0 - \Ba \cdot \nabla \Ncal_\Bp)(\Phi (\By))= O(1)$ as $|\By| \to \infty$. It then follows that
$$
r_0^\dag(\By) = \frac 1 {|\By - \Bp|} \Ba \cdot \nabla \Ncal_\Bp (\Phi(\By)) + O(|\By|^{-1}).
$$
A straightforward computation yields
$$
\frac 1 {|\By - \Bp|} \Ba \cdot \nabla \Ncal_\Bp (\Phi(\By)) = \frac{1} {4\pi} \Ba \cdot (\By-\Bp).
$$ Thus we have
$$
r_0^\dag(\By) = \frac{1} {4\pi}\Ba \cdot (\By-\Bp)  + O(|\By|^{-1}),
$$
or equivalently
\beq\label{condition_r_2}
r_0^* (\Bz) = \frac {\Ge_*} {4\pi} \Ba \cdot (\Bz-\Bp)+ O\left(|\Bz|^{-1}\right) \quad\mbox{as } |\Bz| \rightarrow \infty.
\eeq

Let $w(\Bz):= \Ge_*^{-1} r_0^*(\Bz)$. Then, $w$ also satisfies \eqnref{equation_r_2} and by the maximum principle,
$$
\left|w(\Bz) - {(4\pi)}^{-1} \Ba \cdot (\Bz-\Bp) \right| \leq \norm{{(4\pi)}^{-1} \Ba \cdot (\Bz-\Bp) }_{L^{\infty}(\p D_1 \cup \p D_2)} \leq 1 
$$ 
for all $\Bz \in \mathbb{R}^3 \setminus \overline{D_1 \cup D_2}$. If we take a bounded domain containing $\overline{(D_{1} \cup  D_{2})}$, letting it be $B_3$, then $|w|$ is bounded by $2$ in $\overline{B_3}$, and  $w =0 $ on $\p D_1 \cup \p D_2$. We may use Theorem 1.1 in \cite{LBLY-QJAM-13} and Lemmas 2.1 and 2.3  in \cite{BLY-ARMA-09} to infer that
\beq\label{nablaw}
|\nabla w (\Bz)|  \lesssim \exp \left(-\frac A { {\sqrt \Ge} + |\Bz  | } \right)
\eeq for all $\Bz \in {B_3} \setminus \overline{(D_{1} \cup  D_{2})}$ and for some constant $A$. Since $\nabla w (\Bz) = (1/4\pi) \Ba + O(|\Bz|^{-2})$, we infer from the maximum principle that \eqnref{nablaw} holds for all $\Bz \in \Rbb^3 \setminus \overline{(D_{1} \cup  D_{2})}$. Thus, we have
\beq\label{expone}
|\nabla r_0^* (\Bz)|  \lesssim \Ge_* \exp \left(-\frac A { {\sqrt \Ge} + |\Bz  | } \right).
\eeq
We emphasize here that the constant $A$ can be chosen independently of $\Bp$.

We now prove
\beq\label{exptwo}
| r_0^* (\Bz)|  \lesssim \Ge_* \exp \left(-\frac A { {\sqrt \Ge} +|\Bp|+ |\Bz  |} \right) \left( |\Bz - \Bp| + |\Bp|^2 + \Ge \right)
\eeq
for all $\Bz\in  \Rbb^3 \setminus \overline{(D_{1} \cup  D_{2})}$ if $\Ge$ is sufficiently small. Here and afterwards, the constant $A$ may differ at each occurrence. In fact, since $D_j$ are balls, one can choose a constant $C >0$ regardless of $\Bp$ with $|\Bp|\leq 1/2$ and small $\Ge>0$ so that for any point $\Bz\in  \Rbb^3 \setminus \overline{(D_{1} \cup  D_{2})}$ there are two points, say $\Bz_1$ and $\Bz_2$ in $B_{C(|\Bz| + \sqrt \Ge +|\Bp|)}$ such that the line segments $\ol{\Bz\Bz_1}$, $\ol{\Bz_1\Bz_2}$, $\ol{\Bz_2\Bp}$ lie in $B_{C(|\Bz| + \sqrt \Ge +|\Bp|)} \setminus \ol{(D_1 \cup D_2)}$, and the following also holds independently of $\Bz$ :
$$
|\Bz-\Bz_1|+|\Bz_1-\Bz_2|+|\Bz_2-\Bp|\le  C |\Bz-\Bp|.
$$ 
In addition, there is a point ${\tilde \Bp} $ on $\p D_1$ such that 
$$
|\Bp -{\tilde \Bp}| = \mbox{dist}(\Bp, \p D_1) \simeq |\Bp|^2 + \Ge .
$$
Since $r_0^*(\tilde \Bp)=0$, we may apply mean value theorem on each line segments including $\ol{\Bp {\tilde\Bp}}$ to obtain \eqnref{exptwo} from \eqnref{expone} with a different constant $A$.

Since $r_0 (\Bx ) =  |\Bx - \Bp|^{-1} r_0^\dag (\Phi(\Bx))$ and $|\nabla \Phi(\Bx)| \lesssim |\Bx-\Bp|^{-2}$, we have
\beq\label{2.11}
|\nabla r_0(\Bx)| \lesssim \frac 1 {|\Bx - \Bp|^3}  \left| ( \nabla r_0^\dag)(\Phi(\Bx))  \right|+ \frac 1 {|\Bx - \Bp|^2} \left|r_0^\dag (\Phi(\Bx))\right|.
\eeq
Since $\Phi(\Bx)=\By= \Ge_* \Bz + \Bp_*$, it follows from \eqnref{dagsharp}, \eqnref{expone} and \eqnref{exptwo} that
\begin{align*}
\left| ( \nabla r_0^\dag)(\Phi(\Bx))  \right| = \frac 1 {\Ge_*} \left| ( \nabla r_0^*)(\Bz))  \right| \lesssim \exp \left(-\frac {A} { {\sqrt \Ge} +|\Bp|+ |\Bz  | } \right),
\end{align*}
and
\begin{align*}
\left|  r_0^\dag(\Phi(\Bx))  \right| &\lesssim \Ge_* (|\Bz-\Bp| + |\Bp|^2 + \Ge) \exp \left(-\frac {A} { {\sqrt \Ge}+|\Bp| + |\Bz  | } \right) \\& \lesssim \left(\frac 1 {|\Bx - \Bp |}  + 1\right)  \exp \left(-\frac {A} { {\sqrt \Ge}+|\Bp| + |\Bz  | } \right).
\end{align*}
Thus \eqnref{2.11} yields
$$
|\nabla r_0(\Bx)| \lesssim \left(\frac 1 {|\Bx - \Bp|^3} + \frac 1 {|\Bx - \Bp|^2}  \right)\exp \left(-\frac {A} { {\sqrt \Ge} +|\Bp|+ |\Bz  | } \right).
$$ This inequality can be rewritten as
\beq
|\nabla r_0(\Bx)| \lesssim\left(\frac 1 {|\Bx - \Bp|^3} + \frac 1 {|\Bx - \Bp|^2}  \right)  \exp \left(-\frac {A} { {\sqrt \Ge}  +|\Bp| + |\Bz-\Bp  | } \right). 
\eeq
We see from \eqnref{xyz} that $\Ge_* \approx (\Ge+|\Bp|^2)^{-1}$ and
$$
|\By - \Bp| = \Ge_* |\Bz - \Bp| \approx (\Ge+|\Bp|^2)^{-1} |\Bz - \Bp|.
$$
Therefore, we have
$$
|\Bz-\Bp| \approx \frac{\Ge+|\Bp|^2}{|\Bx-\Bp|},
$$
and \eqnref{case12:gamma} follows. This completes the proof. \qed

The following lemma is used to derive estimates \eqref {farther} and \eqref{est0012} for $|\nabla u (\Bx)|$ for $\Bx$ with $|\Bx| >4$. Even if it may be well-known, we include a short proof based on Uns{\"o}ld's theorem to clarify the dependence of the constant.

\begin{lem} \label{lem_tau}
Let $w$ be a harmonic function defined in $\Rbb^3 \setminus  B_{3}$. If  $|w(\Bx)| = O(|\Bx|^{-2})$ as $|\Bx| \rightarrow \infty$, then, there is a constant $C$ such that
$$
|\nabla w (\Bx)| \leq C \norm {\nabla w }_{L^{2}(\p B_3)}  |\Bx|^{-3}
$$
for all $\Bx$ with $|\Bx| > 4$.
\end{lem}

\pf  
It is well-known that $\nabla w(\Bx)= O(|\Bx|^{-3}$. For $j=1,2,3$, a partial derivative $ \p_j w $ can be written in terms of spherical harmonics as
$$
\p_j w(\Bx) = \sum_{l= 2}^{\infty} \frac 1 {|\Bx|^{l+1}}\sum_{m= -l}^{l}c_{lm} Y_l ^m (\Gt, \Gvf)
$$
for all $\Bx$ with $|\Bx| \geq 3$ (the coefficient depends on $j$). Here, $\Bx = |\Bx|\left(  \cos\Gt \sin\Gvf ,\sin\Gt \sin\Gvf,\cos\Gvf\right)$.   Since $|Y_l ^m(\Gt, \Gvf)| \leq \frac {4\pi} {2l +1}$ by Uns{\"o}ld's theorem, we have
\begin{align*} 
| \p_j w(\Bx) |=& \left( \frac 3  {\left| \Bx\right|}\right)^{3}  \left|  \sum_{l= 2}^{\infty} \left( \frac 3  {\left| \Bx\right|}\right)^{l-2}  \frac 1 {3 ^{l+1}} \sum_{m= -l}^{l}c_{lm} Y_l ^m (\Gt, \Gvf) \right| \\& \lesssim  \frac 1  {\left| \Bx\right|^{3}} \sqrt {\sum_{l=2}^{\infty} \sum_{m=-l}^l\left( \frac 1 {3^{l+1}}  |c_{lm}|\right)^2} \sqrt {\sum_{l=2}^{\infty} \sum_{m=-l}^l\left( \left(\frac 3 4\right)^{l-2}   \frac {4\pi} {2l +1} \right)^2}   \\& \lesssim  \frac 1  {\left| \Bx\right|^{3}} \norm{\p_j w}_{L^{2} (\p B_{3})}   
\end{align*}
for all $\Bx$ with $|\Bx| \geq 4$.   This completes the proof.  \qed

\section{Proofs of Theorems \ref{thm100} and \ref{thm010}}

In this section we deal with the cases when $\Ba=(1,0,0)$ or $(0,1,0)$. As we have seen in the previous section, the function $r$, which is the non-enhanced part in the decomposition \eqnref{r_def} and the solution to \eqnref{eqn_r}, attains $0$ on both $\p D_1$ and $\p D_2$. We first derive estimates for $r$.

\subsection{Estimates for r}

We prove the following proposition.

\begin{prop}\label{prop:case12}
If $\Ba$ is either $(1,0,0)$ or $(0,1,0)$, then the following estimate holds for $\nabla r$, where $r$ is the solution to \eqnref{r_def}:
There exist positive constants $A$ and $\Ge_0$ such that
\beq\label{case12:r}
|\nabla r (\Bx)|  \lesssim \frac 1 {|\Bx - \Bp|^3} \exp \left(- \frac {A|\Bx- \Bp|}{( \sqrt \Ge + |\Bp|)  |\Bx - \Bp | + \Ge + |\Bp|^2 }\right)
\eeq
for all $\Bx \in \Rbb^3 \setminus \overline{D_1 \cup D_2}$, all $|\Bp| \le 1/2$, and all $\Ge \le \Ge_0$.
\end{prop}

The estimate \eqnref{case12:r} already has an important implication. It show that the solution $r$ to \eqnref{eqn_r} when  $c=0$ and $|\Bp| \lesssim \sqrt \Ge$ decays exponentially fast away from the location $\Bp$ of the emitter as the distance $\Ge$ between two inclusions tends to $0$ (see remark after the statement of Theorem \ref{thm010} in Introduction).

\medskip

\noindent{\sl Proof of Proposition \ref{prop:case12}}.
Since $r=0$ on $\p D_1 \cup \p D_2$, $r=r_0$ in $\Rbb^3 \setminus (D_1 \cup D_2)$ where $r_0$ is the solution to \eqnref{eqn_r_0_def}. Lemma \ref{lem:case123} implies that  \eqref{case12:r} holds in a bounded region, say $B_4 \setminus (D_1 \cup D_2)$. One can also see from the same lemma that for all $\Bx \in \p B_3$
\beq
\norm {\nabla r }_{L^{\infty}(\p B_3)} \lesssim \exp\left(- \frac A{2(\sqrt \Ge + |\Bp|) }\right) .
\eeq
Since $r(\Bx) = O(|\Bx|^{-2})$ as $|\Bx| \rightarrow \infty$, Lemma \ref{lem_tau} yields that
\begin{align*}
|\nabla r (\Bx)|&\lesssim \frac 1 {|\Bx|^{3}}  \exp\left(- \frac A{2(\sqrt \Ge + |\Bp|) }\right) \\
&\lesssim \frac 1 {|\Bx - \Bp|^3} \exp \left(- \frac {A|\Bx- \Bp|}{(\sqrt \Ge + |\Bp|)  |\Bx - \Bp | + \Ge + |\Bp|^2 }\right)
\end{align*}
for all $\Bx \in \Rbb^3 \setminus B_4$, where the constant $A$ may differ at each occurrence. Thus, \eqref{case12:r} holds in $ \Rbb^3  \setminus (D_1 \cup D_2)$.  \qed

We are now ready to prove Theorem \ref{thm010}.

\noindent{\sl Proof of Theorem \ref{thm010}}.
As shown in \eqref {symm010},  $Q \equiv 0$ if $\Ba=(0,1,0)$. The global upper bound \eqref {est010} of $|\nabla u|$ in Theorem \ref{thm010} is an immediate consequence of the proposition above.

\par  To show the local estimate \eqref{est010_2} near $\Bp$, we write
$$
u(\Bx)= r(\Bx) =\p_{2} \Ncal_{\Bp}(\Bx) +  (r(\Bx) - \p_{2} \Ncal_{\Bp}(\Bx)).
$$
From definition, one can see
\beq
|\nabla \p_{2} \Ncal_{\Bp}(\Bx)| \simeq |\Bx-\Bp|^{-3}. \label{est010_2_partial_x_2_N}
\eeq
In what follows, we employ the maximum principle to find an upper bound of $ |\nabla (r - \p_{2} \Ncal_{\Bp})|$ near $\Bp$. Choose a constant $M>0$ regardless of $\Ge$ so that
\beq\label{constM}
\overline {B_{M(\Ge+|\Bp|^2)} (\Bp)}\subset \Rbb^3 \setminus \overline {D_1 \cup D_2}
\eeq
for all $\Bp$ with $|\Bp| < 1/2$. As observed in \eqnref{518}, $ \nabla (r - \p_{2} \Ncal_{\Bp})$ is harmonic in $B_{M(\Ge+|\Bp|^2)} (\Bp)$. If $\Bx \in \p B_{M(\Ge+|\Bp|^2)} (\Bp)$, then it follows from Proposition \ref{prop:case12} that
$$
\left|\nabla \left( r(\Bx) - \p_{2} \Ncal_{\Bp}(\Bx)\right)\right|\leq   \left|\nabla r(\Bx)  \right| +  \left|\nabla  \p_{2} \Ncal_{\Bp}(\Bx) \right| \leq C_1 (\Ge+|\Bp|^2)^{-3}
$$
for some $C_1 >0$. Thus by the maximum principle
$$
\left|\nabla \left( r(\Bx) - \p_{2} \Ncal_{\Bp}(\Bx)\right)\right| \leq C_1 (\Ge+|\Bp|^2)^{-3}
$$
for all $\Bx \in  B_{M(\Ge+|\Bp|^2)} (\Bp)$.  This together with \eqref{est010_2_partial_x_2_N} yields Theorem \ref{thm010} (i).
\qed

\subsection{Estimates for Q} \label{Estimates_for_Q}

In this section, we estimate $\nabla Q$ when $\Ba = (1,0,0)$. Theorem \ref{thm100} is proved as a consequence.

First, we have the following lemma.

\begin{lem}\label{lem:Estimates_for_Q}
Let $u$ be the solution to \eqnref{main_equation} with $\Ba = (1,0,0)$. It holds that
\beq\label{pogap}
\left| u|_{\p D_2} - u|_{\p D_1} \right| \simeq \frac 1 {(\Ge+|\Bp|^2) |\log \Ge|},
\eeq
if $|\Bp| < 1$ and $\Ge$ is sufficiently small.
\end{lem}

\pf
Let $q$ be the function defined by \eqnref{qdef}. Thanks to the formula \eqnref{difference}, we have
\beq
u|_{\p D_2} - u|_{\p D_1} = \int_{\p(D_1 \cup D_2)} (\Ba \cdot \nabla \Ncal_\Bp) \p_\nu q \, d\Gs.
\eeq
According to \eqnref{h:sameradii},
$$
\GD q(\Bx)=\frac{1}{\sum_{n=0}^{\infty} q_n} \sum_{n=0}^{\infty} q_n \left(- \Gd_{\Bp_n} + \Gd_{-\Bp_n} \right).
$$
We then apply Green's formula in $D_1 \cup D_2$ to have
$$
u|_{\p D_2} - u|_{\p D_1} = - \int_{D_1 \cup D_2} (\Ba \cdot \nabla \Ncal_\Bp) \GD q \, d\Bx = - \frac 1 {4 \pi} \int_{D_1 \cup D_2} \frac {x_1 }{|\Bx - \Bp|^3} \GD q \, d\Bx,
$$
and hence
\beq\label{3600}
u|_{\p D_2} - u|_{\p D_1} =  \frac 1 {2\pi \sum_{n= 0 }^{\infty}  q_n }  \sum_{n=0} ^{\infty} \frac {q_n p_n}{ (|\Bp|^2 + p_n^2)^{3/2}}.
\eeq

We now show that
\beq\label{4000}
\sum_{n=0} ^{\infty} \frac {q_n p_n}{ (|\Bp|^2 + p_n^2)^{3/2}} \simeq \frac 1 { \Ge + |\Bp|^2}.
\eeq
In fact, letting $C >1$ be the constant appeared in \eqnref{3700}, we have 
$$
\sum_{n=0} ^{\infty} \frac {q_n p_n}{ (|\Bp|^2 + p_n^2)^{3/2}} \ge \sum_{n \leq 1/(C\sqrt{\Ge})} \frac {q_n p_n}{ (|\Bp|^2 + p_n^2)^{3/2}}  \simeq \sum_{n \leq 1/(C\sqrt{\Ge})} \frac {(n+1)^{-2}}{ |\Bp|^3+ (n+1)^{-3}} \simeq \frac 1 { \Ge + |\Bp|^2}
$$ for all $\Bp$ with $|\Bp|< 1$. We also have
\begin{align*}
 \sum_{n=0} ^{\infty} \frac {q_n p_n}{ (|\Bp|^2 + p_n^2)^{3/2}} &=  \sum_{n < 1/(C\sqrt{\Ge})} \frac {q_n p_n}{ (|\Bp|^2 + p_n^2)^{3/2}} +  \sum_{n \geq 1/{C\sqrt{\Ge}}} \frac {q_n p_n}{ (|\Bp|^2 + p_n^2)^{3/2}}  \\
 & \lesssim  \frac 1 { \Ge + |\Bp|^2} + \left( \frac {\sqrt {\Ge}} { (\Ge+ |\Bp|^2 )^{3/2}} \sum_{n \geq 1/(C\sqrt{\Ge})} q_n \right) \lesssim  \frac 1 { \Ge + |\Bp|^2} ,\end{align*}
where the last inequality follows from \eqnref{3900}. Thus \eqnref{4000} follows.

The desired estimate \eqnref{pogap} is an immediate consequence of \eqref{eq:sum_q_n}, \eqnref{3600} and \eqnref{4000}.
\qed

\medskip
\noindent{\sl Proof of Theorem \ref{thm100}}.
We may take $1/4$ for $C_0$ in Theorem \ref{thm100} so that
\beq\label{1/4}
|\Bp| \leq (1/4) |\log \Ge|^{-2}.
\eeq

If $|\Bx| \leq 2 |\log \Ge|^{-2}$ and $\Bx \in \Rbb^3 \setminus \overline{D_1\cup D_2}$, it was shown in \cite{KLY-SIAP-14} that
\beq\label{Q_sub_ineq_1}
\frac{1}{2 q|_{\p D_2}} |\nabla q (\Bx)| \simeq \frac 1 {\Ge + x_2^2 + x_3 ^2} \simeq \frac 1 {\Ge + |\Bx|^2} .
\eeq
Therefore, we can infer from the definition \eqnref{Qdef} of $Q$ and \eqnref{pogap} that
\beq\label{Q_sub_ineq_ok}
|\nabla Q (\Bx)| \simeq \frac 1 {(\Ge + |\Bp|^2) |\log \Ge| (\Ge + |\Bx|^2)}.
\eeq

If $|\Bx-\Bp| \leq  C  (\Ge + |\Bp|^2)$ for a constant $C>0$ and $\Bx \in \Rbb^3 \setminus \overline{(D_{1} \cup  D_{2})}$, then
\beq\label{GeBx}
\Ge + |\Bx|^2 \simeq \Ge + |\Bp|^2.
\eeq
In this case, we write
\beq\label{111}
u(\Bx)= Q(\Bx)+ r(\Bx) = Q(\Bx)+ (r(\Bx) - \p_{1} \Ncal_{\Bp}(\Bx)) + \p_{1} \Ncal_{\Bp}(\Bx),
\eeq
and estimate three terms on the right-hand side one by one.
By explicit computations, one can see
\beq\label{112}
|\nabla \p_{1} \Ncal_{\Bp}(\Bx)| \simeq |\Bx-\Bp|^{-3}.
\eeq
We also see from \eqnref{Q_sub_ineq_ok} and \eqnref{GeBx} that
\beq\label{113}
|\nabla Q (\Bx)| \lesssim \frac 1 {(\Ge + |\Bp|^2)^2  |\log \Ge|} \lesssim \frac{\Ge+ |\Bp|^2}{|\log \Ge|}|\Bx-\Bp|^{-3}.
\eeq
Let $M$ be the constant appearing in \eqnref{constM} and suppose that $C< M$.  Since $|\Bx-\Bp| \leq  C  (\Ge + |\Bp|^2)$,
we see from the maximum principle and  \eqnref{case12:r} that
\begin{align*}
|\nabla (r - \p_{1} \Ncal_{\Bp})(\Bx)| &\leq \norm{\nabla r}_{L^{\infty}(\p B_{M (\Ge + |\Bp|^2) } (\Bp )) }  +  \norm{\nabla \p_{1} \Ncal_{\Bp} }_{L^{\infty} (\p B_{M (\Ge + |\Bp|^2) } (\Bp ))}\\
&  \lesssim \frac 1 {|\Bx-\Bp|^{3}}. 
\end{align*}
Thus (i) hold for all $C<M$.

If $C (\Ge +|\Bp|^2) \leq |\Bx-\Bp| \leq C_1 (\Ge +|\Bp|^2 )|\log \Ge|$ for some $C_1>C$ and $\Bx \in \Rbb^3 \setminus \overline{(D_{1} \cup  D_{2})}$, then
$$
(\sqrt \Ge+|\Bp| ) |\Bx - \Bp | + \Ge+|\Bp| ^2 \lesssim  \Ge+|\Bp| ^2.
$$
Thus \eqnref{case12:r} yields
$$
|\nabla r (\Bx)|  \lesssim \frac 1 {|\Bx - \Bp|^3} \exp \left(- A \frac{|\Bx- \Bp|}{\Ge + |\Bp|^2}\right).
$$
(The constant $A$ may differ at each occurrence.) So, \eqnref{between} follows from \eqref{Q_sub_ineq_ok}, since $\Ge + |\Bp|^2 \simeq \Ge + |\Bx|^2$.

We now show that there exists a constant $C_2 >0$ such that if $C_2(\Ge+|\Bp|^2) |\log \Ge| \leq |\Bx-\Bp| \leq |\log \Ge|^{-2}$, then
\beq\label{twice}
|\nabla r (\Bx)| \le \frac{1}{2} |\nabla Q (\Bx)|,
\eeq
provided that $\Ge$ is sufficiently small. Once this fact is proved, it follows that $|\nabla u (\Bx)| \simeq |\nabla Q(\Bx)|$, and we have \eqnref{far}. We then choose $C$ and $C_1$ in the above so that $C^{-1}=C_1 \ge C_2$. Then we have (i), (ii) and (iii).

To prove \eqnref{twice}, we first consider the case when $3(\sqrt{\Ge} + |\Bp|) \le |\Bx-\Bp| \leq  |\log \Ge|^{-2}$. In this case, we have
$$
|\Bx| \le |\Bx-\Bp| + |\Bp| \lesssim |\Bx-\Bp|,
$$
and hence
$$
|\Bx| + \sqrt{\Ge} + |\Bp| \lesssim |\Bx-\Bp| .
$$
It then follows from \eqnref{case12:r} that
$$
|\nabla r (\Bx)|  \lesssim \frac 1 {|\Bx - \Bp|^3} \lesssim \frac{\sqrt{\Ge + |\Bp|^2} |\log \Ge|}{(\Ge + |\Bp|^2)|\log \Ge|(\Ge + |\Bx|^2)}.
$$
Thanks to \eqnref{1/4}, we have $\sqrt{\Ge + |\Bp|^2} |\log \Ge| \lesssim 1$. Thus \eqnref{twice} follows from \eqref{Q_sub_ineq_ok} in this case.

Secondly, we deal with the case when
\beq\label{upperlower}
C_2(\Ge + |\Bp|^2) |\log \Ge| \leq |\Bx-\Bp| \le 3(\sqrt{\Ge}+|\Bp|),
\eeq
and determine the constant $C_2$. In this case, we have
$$
(\sqrt \Ge + |\Bp|)  |\Bx - \Bp | + \Ge + |\Bp|^2 \le 7(\Ge + |\Bp|^2).
$$
Thus \eqnref{case12:r} yields
\begin{align*}
|\nabla r (\Bx)| &\lesssim  \frac 1 {|\Bx - \Bp |^3} \exp \left(-\frac{A|\Bx-\Bp|}{7(\Ge + |\Bp|^2)}\right) \\
&\lesssim \frac 1 {(\Ge+|\Bp|^2)^3 |\log \Ge|^3} \exp \left( \frac{AC_2}{7} \log \Ge \right) \lesssim \frac { \Ge^{AC_2/7} } {(\Ge+|\Bp|^2)^{3}|\log \Ge|^3}.
\end{align*}
In view of \eqnref{Q_sub_ineq_ok}, we write the above as
$$
|\nabla r (\Bx)| \lesssim  \frac 1 {(\Ge + |\Bp|^2) |\log \Ge| (\Ge + |\Bx|^2)} \frac{\Ge^{AC_2/7} (\Ge + |\Bx|^2)}{(\Ge + |\Bp|^2)^2 |\log \Ge|^2} \simeq |\nabla Q (\Bx)| \frac{\Ge^{AC_2/7} (\Ge + |\Bx|^2)}{(\Ge + |\Bp|^2)^2 |\log \Ge|^2} .
$$
Observe that $\Bx$ in the range given by \eqnref{upperlower} satisfies $|\Bx| \lesssim |\log\Ge|^{-2}$. Thus we have
$$
|\nabla r (\Bx)| \lesssim |\nabla Q (\Bx)| \frac{\Ge^{AC_2/7}}{\Ge^2 |\log \Ge|^{  6}} .
$$
So, if $C_2$ is large enough, for example, if $AC_2/7> 2$, then $\frac{\Ge^{AC_2/7}}{\Ge^2 |\log \Ge|^6}$ can be arbitrarily small (provided that $\Ge$ is small enough), and \eqnref{twice} holds.

Finally, we prove (iv) for the case when $|\Bx - \Bp| \geq |\log \Ge|^{-2}$. We consider two cases separately: when $|\Bx|< 4$ and $|\Bx - \Bp| \geq |\log \Ge|^{-2}$, and when  $|\Bx| \geq 4$.  

In the first case, $|\Bx|\geq (3/4) |\log \Ge|^{-2}$ thanks to \eqnref{1/4}.  Since $Q$ is constant on $\p D_i $, $i=1,2$, one can use the inversion with respect to either $\p D_1$ or $\p D_2$ and the Kelvin transform to see that there is a small constant $C$ independent of $\Ge$ such that for any $\Bx \in B_4 \setminus (D_1 \cup D_2)$,  $Q$ can be locally extended into $B_{C(|\Bx|^2 + \Ge)}(\Bx)$ as a harmonic functions. We may assume $C < 1/20$ by taking even smaller $C$. By the maximum principle, $|Q| \leq Q|_{\p D_2} - Q|_{\p D_1}$ in $\Rbb^3 \setminus (D_1 \cup D_2)$. Since $C < 20^{-1}$, we have for any $\Bx' \in B_{C(|\Bx|^2 + \Ge)} (\Bx)$
$$
|\Bx' - \Bc_j| \geq |\Bx - \Bc_j| - |\Bx' - \Bx| \geq 1 - C (|\Bx|^2 + \Ge) \geq 1- C (16+\Ge) > C,
$$
where $\Bc_j$ is the center of $D_j$. Now we have an estimate for the Kelvin transform of $Q$: 
$$
\left| \frac 1 {|\Bx'- \Bc_j| } Q \left( \frac {\Bx'- \Bc_i}{|\Bx'- \Bc_i|^2}+ \Bc_i\right)\right|\leq \frac 1 C (Q|_{\p D_2} - Q|_{\p D_1})
$$ 
for $j=1,2$, and hence, the extended function is bounded by $C^{-1} (Q|_{\p D_2} - Q|_{\p D_1})$. We then infer from the standard gradient estimate for harmonic functions that
$$
|\nabla Q (\Bx)| \lesssim \frac {Q|_{\p D_2} - Q|_{\p D_1} } {C^{2} (|\Bx|^2 + \Ge) } .
$$
Since $Q|_{\p D_2} - Q|_{\p D_1} = u|_{\p D_2} - u|_{\p D_1}$, Lemma \ref{lem:Estimates_for_Q} yields
\beq\label{estQ}
|\nabla Q (\Bx)| \lesssim \frac {1}{ (\Ge + |\Bp|^2) |\log \Ge||\Bx|^2}
\eeq
for all $\Bx \in B_4 \setminus (D_1 \cup D_2)$ with  $|\Bx - \Bp| \geq |\log \Ge|^{-2}$.

Since $| \Bx - \Bp | \geq |\log \Ge|^{-2}$, we have $|\Bx| \geq (3/4) |\log \Ge|^{-2} \geq 3 |\Bp|$, $| \Bx - \Bp | \geq 2\Ge |\log \Ge|^{3/2}$, and $|\Bx-\Bp| \geq 4  |\Bp| \geq 2 |\Bp|^2 |\log \Ge|^{3/2} $, provided that $\Ge$ is small enough. Thus,
$$
| \Bx - \Bp | \geq  | \Bx| - |\Bp | \geq  (2/3) |\Bx| \quad\mbox{and}\quad |\Bx - \Bp| \geq (\Ge + |\Bp|^2) |\log  \Ge|^{3/2}.
$$
By Proposition \ref{prop:case12},
$$
| \nabla r (\Bx)| \lesssim \frac 1 {|\Bx - \Bp|^3} \lesssim \frac {{|\log \Ge|^{-1/2}}}{ (\Ge + |\Bp|^2)|\log \Ge||\Bx|^2}.
$$
This together with \eqnref{estQ} yields 
\beq\label{farther_100_1st_case}
| \nabla u (\Bx)  | \leq  |\nabla Q (\Bx)| + |\nabla r(\Bx) | \lesssim  \frac {1}{ (\Ge + |\Bp|^2) |\log \Ge||\Bx|^2}
\eeq
for all $\Bx$ with $ |\Bx|< 4$ and  $|\Bx - \Bp| \geq |\log \Ge|^{-2}$.

Suppose now that $|\Bx|\geq 4$. By \eqnref{farther_100_1st_case}, we have
$$
\norm{\nabla u}_{L^{2} (\p B_{3})} \lesssim \frac {1}{ (\Ge + |\Bp|^2) |\log \Ge|} .$$  By Lemma \ref{lem_tau},
$$
|\nabla u(\Bx) | \lesssim  \frac 1  {\left| \Bx\right|^{3}} \norm{\nabla  u}_{L^{2} (\p B_{3})} \lesssim  \frac {1}{ (\Ge + |\Bp|^2) |\log \Ge| |\Bx|^3}
$$
for all $\Bx$ with $|\Bx| \geq 4$. The proof is complete.
\qed

\section{Proof of Theorem \ref{thm001}}\label{sec001}

In this section we deal with the case when $\Ba = (0,0,1)$ to prove Theorem \ref{thm001}. By (iii) in section \ref{sec:decom}, $Q \equiv 0$ and $u|_{\p D_1}=u|_{\p D_2}$. We begin with the following lemma.

\begin{lem}
Let $u$ be the solution to \eqnref{main_equation} with $\Ba=(0,0,1)$. If $|\Bp|< 1/2$, then
\beq\label{boundeu}
|u|_{\p D_1} | \lesssim  1.
\eeq
\end{lem}

\pf
Since $u|_{\p D_1}=u|_{\p D_2}$, it follows from \eqnref{sum} that
$$
u|_{\p D_1} = \frac{1}{2} \int_{\p(D_1 \cup D_2)} (\p_{3} \Ncal_{\Bp}) \p_\nu q^\perp  \, d\Gs.
$$
Then, the relation \eqnref{qperp} yields
$$
u|_{\p D_1} = \frac 1 2 \int_{\p(D_1 \cup D_2)} (\p_{3} \Ncal_{\Bp}) \p_\nu (v_f - 2 \Ncal_{\Bc_1} + q)  \, d\Gs .
$$
Note that
\beq\label{zeroint}
\int_{\p(D_1 \cup D_2)} (\p_{3} \Ncal_{\Bp}) \p_\nu q \, d\Gs = 0.
\eeq
In fact, since $q$ is skew-symmetric with respect to the $x_2x_3$-plane, so is $\p_\nu q$. On the other hand, $\p_{3} \Ncal_{\Bp}$ is symmetric with respect to the $x_2x_3$-plane. So, \eqnref{zeroint} follows. Thus
\beq\label{uoneform}
u|_{\p D_1}=  \frac 1 2 \int_{\p(D_1 \cup D_2)} (\p_{3} \Ncal_{\Bp}) \p_\nu (v_f - 2 \Ncal_{\Bc_1})  \, d\Gs.
\eeq

We claim that $v_f$ is bounded  regardless of $\Ge$. In fact, since $v_f(\Bx)=O(|\Bx|^{-2})$ as $|\Bx| \to \infty$, $v_f$ attains its maximum and minimum on $\p (D_1 \cup D_2)$. Moreover, the maximum is positive and the minimum is negative. Since $v_f - 2 \Ncal_{\Bc_1}$ is constant, say $\Gl_j$, on $\p D_j$, we have
$$
\max_{\Bx \in \p (D_1 \cup D_2)} |v_f(\Bx)| \le 4 \max_{\Bx \in \p (D_1 \cup D_2)} |\Ncal_{\Bc_1}(\Bx)| + |\Gl_2-\Gl_1| \le C + |\Gl_2-\Gl_1|
$$
for some constant $C$. By \eqnref{difference}, we have
$$
\Gl_2 - \Gl_1  = - 2\int_{\p(D_1 \cup D_2)}  \Ncal_{\Bc_1}(\Bx) \p_\nu q(\Bx)  \, d\Gs ,
$$
and hence
$$
|\Gl_2 - \Gl_1| \leq  2 \int_{\p(D_1 \cup D_2)} | \Ncal_{\Bc_1}(\Bx) \p_\nu q(\Bx)| \, d\Gs
\lesssim \int_{\p(D_1 \cup D_2)}  |\p_\nu q(\Bx)|  \, d\Gs .
$$
Since the constant values of $q$ on $\p D_1$ and $\p D_2$ are the minimum and maximum of $q$, $\p_\nu q$ is either strictly positive or negative on $\p D_j$ by Hopf's lemma. Thus,
$$
|\Gl_2 - \Gl_1|
\lesssim \sum_{j=1}^2 \left| \int_{\p D_j}  \p_\nu q(\Bx)  \, d\Gs \right| = 2.
$$
Therefore, $v_f$ is bounded.

Define $g$ by
\beq\label{zeroint3}
v_f - 2 \Ncal_{\Bc_1} = g + \frac{\Gl_2-\Gl_1}{2 q|_{\p D_2}} q
\eeq
in $\Rbb^3 \setminus (D_1 \cup D_2)$. Then
$$
|g| \le |v_f - 2 \Ncal_{\Bc_1}| + \left| \frac{\Gl_2-\Gl_1}{2 q|_{\p D_2}} q \right| \lesssim 1 + | \Gl_2-\Gl_1| \lesssim 1.
$$
We then obtain in the same way as to derive \eqnref{nablaw} that 
\beq\label{lemma_case3_grad_h}
\left| \nabla g (\Bx) \right| \lesssim \exp \left(-  \frac {A}{\sqrt {\Ge} + |\Bx|}\right)
\eeq
for some constant $A >0$.  

Thanks to \eqnref{zeroint}, \eqnref{uoneform} and \eqnref{zeroint3}, we have
$$
u|_{\p D_1} = \frac 1 2 \int_{\p(D_1 \cup D_2)} (\p_{x_3} \Ncal_{\Bp}) \p_\nu g  \, d\Gs.
$$
Write
$$
u|_{\p D_1}=I_1+I_2:= \frac 1 2  \int_{(\p D_1 \cup \p D_2)\cap B_{3(\sqrt \Ge+|\Bp|)} }+ \frac 1 2  \int_{(\p D_1 \cup \p D_2) \setminus B_{3(\sqrt \Ge+|\Bp|)}}  (\p_{x_3} \Ncal_{\Bp}) \p_{\nu}  g d\Gs.
$$
Then using \eqnref{lemma_case3_grad_h}, we obtain
\begin{align*}
|I_1| &\lesssim \int_{(\p D_1 \cup \p D_2)\cap B_{3(\sqrt \Ge+|\Bp|)}}  \frac 1 {(\Ge + |\Bp|^2)^2} \exp \left(- \frac A {4(\sqrt \Ge + |\Bp|)}\right) d\Gs\\&
\lesssim \frac 1 {\sqrt \Ge + |\Bp|} \exp \left(- \frac A {4(\sqrt \Ge + |\Bp|)}\right) \lesssim 1
\end{align*}
and
\begin{align*}
|I_2| \lesssim \int_{(\p D_1 \cup \p D_2) \setminus B_{3(\sqrt \Ge+|\Bp|)}}  \frac 1 {|\Bx|^2} \exp \left(- \frac A {2|\Bx|}\right) d\Gs
&\lesssim \int_{\sqrt \Ge + |\Bp|} ^{1} \frac 1 {{\rho}^2} \exp \left(- \frac A {2\rho} \right) d\rho  \lesssim 1,
\end{align*}
which yield \eqnref{boundeu}. 
This completes the proof.
\qed

\medskip
\medskip
\noindent{\sl Proof of Theorem \ref{thm001}}.
Now, we prove \eqref {est001}. Note that $u=r$ since $Q=0$. We decompose $u$ as
$$
u = r_0+w_0,
$$
where $r_0$ is the solution to \eqref {eqn_r_0_def} and $w_0$ is the solution to
\beq
\begin{cases}
\GD w_0 = 0  \quad&\mbox{ in } \Rbb^3 \setminus \overline {(D_{1} \cup D_{2})}, \\
\ds  w_0 = u |_{\p D_1} (= u |_{\p D_2})  \quad& \mbox { on } \p D_1 \cup \p D_2,\\
\ds  w_0 (\Bx)= O\left(|\Bx|^{-1}\right)  &\mbox{ as } |\Bx| \rightarrow \infty.
\end{cases}
\eeq
Thus, $r_0 (\Bx) + w_0(\Bx) = O(|\Bx|^{-2})$ as $|\Bx| \rightarrow \infty$. By Lemma \ref{lem:case123}, there exists a positive constant $A$ independent of $\Bp$ such that
\beq\label{Case3:eq1}
|\nabla r_0 (\Bx)|  \lesssim \frac 1 {|\Bx - \Bp|^3} \exp \left(- \frac {A|\Bx- \Bp|}{ (\sqrt \Ge + \Bp) |\Bx - \Bp | + \Ge+ |\Bp|^2}\right)
\eeq
for all $\Bx \in B_{4} \setminus \overline{D_1 \cup D_2  \cup \{\Bp\}}$ and all sufficiently small $\Ge$. Since $w_0|_{\p D_1} = w_0|_{\p D_2}$ and $|w_0 (\Bx)|  \lesssim \left|u|_{\p D_1} \right| \lesssim 1$ in $\Rbb^3 \setminus \overline{D_1 \cup D_2}$ thanks to \eqref{boundeu}, one can see in the same way as \eqnref{nablaw} that there is a positive constant $A$ such that
\beq\label{Case3:eq2}
|\nabla w_0 (\Bx)| \lesssim  \exp \left(\frac {-A}{ \sqrt \Ge + |\Bx|}\right)
\eeq
in $\Rbb^3 \setminus \overline{D_1 \cup D_2}$. Therefore, \eqref {est001} follows.

Since $u$ is harmonic in $\Rbb^3 \setminus \ol{D_1 \cup D_2}$ and $u(\Bx)=O(|\Bx|^{-2})$ as $|\Bx| \to \infty$, we have $\nabla u(\Bx)=O(|\Bx|^{-3})$ as $|\Bx| \to \infty$. Since $\norm {\nabla u}_{L^{2}(\p  B_3)} \lesssim 1$ by \eqref {est001}, Lemma \ref{lem_tau} yields \eqnref{est0012}.

The estimate \eqref{est0013} can be proved in the same as \eqref{est010_2}.
\qed

 \end{document}